\expandafter\chardef\csname pre amssym.def at\endcsname=\the\catcode`\@
\catcode`\@=11

\def\undefine#1{\let#1\undefined}
\def\newsymbol#1#2#3#4#5{\let\next@\relax
 \ifnum#2=\@ne\let\next@\msafam@\else
 \ifnum#2=\tw@\let\next@\msbfam@\fi\fi
 \mathchardef#1="#3\next@#4#5}
\def\mathhexbox@#1#2#3{\relax
 \ifmmode\mathpalette{}{\m@th\mathchar"#1#2#3}%
 \else\leavevmode\hbox{$\m@th\mathchar"#1#2#3$}\fi}
\def\hexnumber@#1{\ifcase#1 0\or 1\or 2\or 3\or 4\or 5\or 6\or 7\or 8\or
 9\or A\or B\or C\or D\or E\or F\fi}

\font\tenmsa=msam10
\font\sevenmsa=msam7
\font\fivemsa=msam5
\newfam\msafam
\textfont\msafam=\tenmsa
\scriptfont\msafam=\sevenmsa
\scriptscriptfont\msafam=\fivemsa
\edef\msafam@{\hexnumber@\msafam}
\mathchardef\dabar@"0\msafam@39
\def\dashrightarrow{\mathrel{\dabar@\dabar@\mathchar"0\msafam@4B}}
\def\dashleftarrow{\mathrel{\mathchar"0\msafam@4C\dabar@\dabar@}}

\def\ulcorner{\delimiter"4\msafam@70\msafam@70 }
\def\urcorner{\delimiter"5\msafam@71\msafam@71 }
\def\llcorner{\delimiter"4\msafam@78\msafam@78 }
\def\lrcorner{\delimiter"5\msafam@79\msafam@79 }
\def\yen{{\mathhexbox@\msafam@55 }}
\def\checkmark{{\mathhexbox@\msafam@58 }}
\def\circledR{{\mathhexbox@\msafam@72 }}
\def\maltese{{\mathhexbox@\msafam@7A }}

\font\tenmsb=msbm10
\font\sevenmsb=msbm7
\font\fivemsb=msbm5
\newfam\msbfam
\textfont\msbfam=\tenmsb
\scriptfont\msbfam=\sevenmsb
\scriptscriptfont\msbfam=\fivemsb
\edef\msbfam@{\hexnumber@\msbfam}

\catcode`\@=\csname pre amssym.def at\endcsname

\expandafter\ifx\csname pre amssym.tex at\endcsname\relax \else \endinput\fi
\expandafter\chardef\csname pre amssym.tex at\endcsname=\the\catcode`\@
\catcode`\@=11
\newsymbol\boxdot 1200
\newsymbol\boxplus 1201
\newsymbol\boxtimes 1202
\newsymbol\square 1003
\newsymbol\blacksquare 1004
\newsymbol\centerdot 1205
\newsymbol\lozenge 1006
\newsymbol\blacklozenge 1007
\newsymbol\circlearrowright 1308
\newsymbol\circlearrowleft 1309
\undefine\rightleftharpoons
\newsymbol\rightleftharpoons 130A
\newsymbol\leftrightharpoons 130B
\newsymbol\boxminus 120C
\newsymbol\Vdash 130D
\newsymbol\Vvdash 130E
\newsymbol\vDash 130F
\newsymbol\twoheadrightarrow 1310
\newsymbol\twoheadleftarrow 1311
\newsymbol\leftleftarrows 1312
\newsymbol\rightrightarrows 1313
\newsymbol\upuparrows 1314
\newsymbol\downdownarrows 1315
\newsymbol\upharpoonright 1316
 
\newsymbol\downharpoonright 1317
\newsymbol\upharpoonleft 1318
\newsymbol\downharpoonleft 1319
\newsymbol\rightarrowtail 131A
\newsymbol\leftarrowtail 131B
\newsymbol\leftrightarrows 131C
\newsymbol\rightleftarrows 131D
\newsymbol\Lsh 131E
\newsymbol\Rsh 131F
\newsymbol\rightsquigarrow 1320
\newsymbol\leftrightsquigarrow 1321
\newsymbol\looparrowleft 1322
\newsymbol\looparrowright 1323
\newsymbol\circeq 1324
\newsymbol\succsim 1325
\newsymbol\gtrsim 1326
\newsymbol\gtrapprox 1327
\newsymbol\multimap 1328
\newsymbol\therefore 1329
\newsymbol\because 132A
\newsymbol\doteqdot 132B
 
\newsymbol\triangleq 132C
\newsymbol\precsim 132D
\newsymbol\lesssim 132E
\newsymbol\lessapprox 132F
\newsymbol\eqslantless 1330
\newsymbol\eqslantgtr 1331
\newsymbol\curlyeqprec 1332
\newsymbol\curlyeqsucc 1333
\newsymbol\preccurlyeq 1334
\newsymbol\leqq 1335
\newsymbol\leqslant 1336
\newsymbol\lessgtr 1337
\newsymbol\backprime 1038
\newsymbol\risingdotseq 133A
\newsymbol\fallingdotseq 133B
\newsymbol\succcurlyeq 133C
\newsymbol\geqq 133D
\newsymbol\geqslant 133E
\newsymbol\gtrless 133F
\newsymbol\sqsubset 1340
\newsymbol\sqsupset 1341
\newsymbol\vartriangleright 1342
\newsymbol\vartriangleleft 1343
\newsymbol\trianglerighteq 1344
\newsymbol\trianglelefteq 1345
\newsymbol\bigstar 1046
\newsymbol\between 1347
\newsymbol\blacktriangledown 1048
\newsymbol\blacktriangleright 1349
\newsymbol\blacktriangleleft 134A
\newsymbol\vartriangle 134D
\newsymbol\blacktriangle 104E
\newsymbol\triangledown 104F
\newsymbol\eqcirc 1350
\newsymbol\lesseqgtr 1351
\newsymbol\gtreqless 1352
\newsymbol\lesseqqgtr 1353
\newsymbol\gtreqqless 1354
\newsymbol\Rrightarrow 1356
\newsymbol\Lleftarrow 1357
\newsymbol\veebar 1259
\newsymbol\barwedge 125A
\newsymbol\doublebarwedge 125B
\undefine\angle
\newsymbol\angle 105C
\newsymbol\measuredangle 105D
\newsymbol\sphericalangle 105E
\newsymbol\varpropto 135F
\newsymbol\smallsmile 1360
\newsymbol\smallfrown 1361
\newsymbol\Subset 1362
\newsymbol\Supset 1363
\newsymbol\Cup 1264
 
\newsymbol\Cap 1265
 
\newsymbol\curlywedge 1266
\newsymbol\curlyvee 1267
\newsymbol\leftthreetimes 1268
\newsymbol\rightthreetimes 1269
\newsymbol\subseteqq 136A
\newsymbol\supseteqq 136B
\newsymbol\bumpeq 136C
\newsymbol\Bumpeq 136D
\newsymbol\lll 136E
 
\newsymbol\ggg 136F
 
\newsymbol\circledS 1073
\newsymbol\pitchfork 1374
\newsymbol\dotplus 1275
\newsymbol\backsim 1376
\newsymbol\backsimeq 1377
\newsymbol\complement 107B
\newsymbol\intercal 127C
\newsymbol\circledcirc 127D
\newsymbol\circledast 127E
\newsymbol\circleddash 127F
\newsymbol\lvertneqq 2300
\newsymbol\gvertneqq 2301
\newsymbol\nleq 2302
\newsymbol\ngeq 2303
\newsymbol\nless 2304
\newsymbol\ngtr 2305
\newsymbol\nprec 2306
\newsymbol\nsucc 2307
\newsymbol\lneqq 2308
\newsymbol\gneqq 2309
\newsymbol\nleqslant 230A
\newsymbol\ngeqslant 230B
\newsymbol\lneq 230C
\newsymbol\gneq 230D
\newsymbol\npreceq 230E
\newsymbol\nsucceq 230F
\newsymbol\precnsim 2310
\newsymbol\succnsim 2311
\newsymbol\lnsim 2312
\newsymbol\gnsim 2313
\newsymbol\nleqq 2314
\newsymbol\ngeqq 2315
\newsymbol\precneqq 2316
\newsymbol\succneqq 2317
\newsymbol\precnapprox 2318
\newsymbol\succnapprox 2319
\newsymbol\lnapprox 231A
\newsymbol\gnapprox 231B
\newsymbol\nsim 231C
\newsymbol\ncong 231D
\newsymbol\diagup 231E
\newsymbol\diagdown 231F
\newsymbol\varsubsetneq 2320
\newsymbol\varsupsetneq 2321
\newsymbol\nsubseteqq 2322
\newsymbol\nsupseteqq 2323
\newsymbol\subsetneqq 2324
\newsymbol\supsetneqq 2325
\newsymbol\varsubsetneqq 2326
\newsymbol\varsupsetneqq 2327
\newsymbol\subsetneq 2328
\newsymbol\supsetneq 2329
\newsymbol\nsubseteq 232A
\newsymbol\nsupseteq 232B
\newsymbol\nparallel 232C
\newsymbol\nmid 232D
\newsymbol\nshortmid 232E
\newsymbol\nshortparallel 232F
\newsymbol\nvdash 2330
\newsymbol\nVdash 2331
\newsymbol\nvDash 2332
\newsymbol\nVDash 2333
\newsymbol\ntrianglerighteq 2334
\newsymbol\ntrianglelefteq 2335
\newsymbol\ntriangleleft 2336
\newsymbol\ntriangleright 2337
\newsymbol\nleftarrow 2338
\newsymbol\nrightarrow 2339
\newsymbol\nLeftarrow 233A
\newsymbol\nRightarrow 233B
\newsymbol\nLeftrightarrow 233C
\newsymbol\nleftrightarrow 233D
\newsymbol\divideontimes 223E
\newsymbol\varnothing 203F
\newsymbol\nexists 2040
\newsymbol\Finv 2060
\newsymbol\Game 2061
\newsymbol\mho 2066
\newsymbol\eth 2067
\newsymbol\eqsim 2368
\newsymbol\beth 2069
\newsymbol\gimel 206A
\newsymbol\daleth 206B
\newsymbol\lessdot 236C
\newsymbol\gtrdot 236D
\newsymbol\ltimes 226E
\newsymbol\rtimes 226F
\newsymbol\shortmid 2370
\newsymbol\shortparallel 2371
\newsymbol\smallsetminus 2272
\newsymbol\thicksim 2373
\newsymbol\thickapprox 2374
\newsymbol\approxeq 2375
\newsymbol\succapprox 2376
\newsymbol\precapprox 2377
\newsymbol\curvearrowleft 2378
\newsymbol\curvearrowright 2379
\newsymbol\digamma 207A
\newsymbol\varkappa 207B
\newsymbol\Bbbk 207C
\newsymbol\hslash 207D
\undefine\hbar
\newsymbol\hbar 207E
\newsymbol\backepsilon 237F
\catcode`\@=\csname pre amssym.tex at\endcsname

\magnification=1200
\hsize=468truept
\vsize=646truept
\voffset=-10pt
\parskip=4pt
\baselineskip=14truept
\count0=1

\dimen100=\hsize

\def\leftill#1#2#3#4{
\medskip
\line{$
\vcenter{
\hsize = #1truept \hrule\hbox{\vrule\hbox to  \hsize{\hss \vbox{\vskip#2truept
\hbox{{\copy100 \the\count105}: #3}\vskip2truept}\hss }
\vrule}\hrule}
\dimen110=\dimen100
\advance\dimen110 by -36truept
\advance\dimen110 by -#1truept
\hss \vcenter{\hsize = \dimen110
\medskip
\noindent { #4\par\medskip}}$}
\advance\count105 by 1
}
\def\rightill#1#2#3#4{
\medskip
\line{
\dimen110=\dimen100
\advance\dimen110 by -36truept
\advance\dimen110 by -#1truept
$\vcenter{\hsize = \dimen110
\medskip
\noindent { #4\par\medskip}}
\hss \vcenter{
\hsize = #1truept \hrule\hbox{\vrule\hbox to  \hsize{\hss \vbox{\vskip#2truept
\hbox{{\copy100 \the\count105}: #3}\vskip2truept}\hss }
\vrule}\hrule}
$}
\advance\count105 by 1
}
\def\midill#1#2#3{\medskip
\line{$\hss
\vcenter{
\hsize = #1truept \hrule\hbox{\vrule\hbox to  \hsize{\hss \vbox{\vskip#2truept
\hbox{{\copy100 \the\count105}: #3}\vskip2truept}\hss }
\vrule}\hrule}
\dimen110=\dimen100
\advance\dimen110 by -36truept
\advance\dimen110 by -#1truept
\hss $}
\advance\count105 by 1
}
\def\insectnum{\copy110\the\count120
\advance\count120 by 1
}

\font\ninerm=cmr9
\font\eightrm=cmr8

\font\tenrm=cmr10 at 10pt

\font\sc=cmcsc10

%\font\sevenmsy=msym7
%\font\fivemsy=msym5
%\newfam\msyfam
%\def\msy{\fam\msyfam\tenmsy}
%\textfont\msyfam=\tenmsy
%\scriptfont\msyfam=\sevenmsy
%\scriptscriptfont\msyfam=\fivemsy
\def\msb{\fam\msbfam\tenmsb}

\def\bbc{{\msb C}}

\def\bbi{{\msb I}}

\def\bbp{{\msb P}}
\def\bbq{{\msb Q}}

\def\bbz{{\msb Z}}

\def\grD{\Delta}

\def\grG{\Gamma}

\def\grL{\Lambda}

\def\gra{\alpha}

\def\gre{\epsilon}

\def\grl{\lambda}

\def\grt{\tau}

\def\grz{\zeta}

\def\la#1{\hbox to #1pc{\leftarrowfill}}
\def\ra#1{\hbox to #1pc{\rightarrowfill}}

\def\fract#1#2{\raise4pt\hbox{$ #1 \atop #2 $}}
\def\decdnar#1{\phantom{\hbox{$\scriptstyle{#1}$}}
\left\downarrow\vbox{\vskip15pt\hbox{$\scriptstyle{#1}$}}\right.}

\def\bowtie{\hbox to 1pt{\hss}\raise.66pt\hbox{$\scriptstyle{>}$}
\kern-4.9pt\triangleleft}
\def\hsmash{\triangleright\kern-4.4pt\raise.66pt\hbox{$\scriptstyle{<}$}}
\def\boxit#1{\vbox{\hrule\hbox{\vrule\kern3pt
\vbox{\kern3pt#1\kern3pt}\kern3pt\vrule}\hrule}}

\def\za{\vrule height6pt width4pt depth1pt}

\font\aa=eufm10

\def\Got#1{\hbox{\aa#1}}

\def\bfw{{\bf w}}

\def\calo{{\cal O}}

\def\calf{{\cal F}}

\def\call{{\cal L}}
\def\calm{{\cal M}}

\def\calo{{\cal O}}

\def\cals{{\cal S}}

\def\calz{{\cal Z}}

\def\gG{{\Got G}}

\font\svtnrm=cmr17

%\font\sbb=msym8
%\font\ninesl=cmsl9
\font\bsc=cmcsc10 at 10truept

\def\coker{\hbox{coker}}
\def\ker{\hbox{ker}}

\def\Ric{\hbox{Ric}}

\def\Se{Sasakian-Einstein }

\def\sea{1}

\def\sec{2}
\def\sed{3}
\def\tvo{4}
\def\see{5}
\settabs 9\columns

\phantom{ooo}
\bigskip\bigskip
\centerline{\svtnrm Sasakian-Einstein Structures on $9\#(S^2\times
S^3)$}  \medskip
%\centerline{$9\#(S^2\times S^3)$}

\bigskip\bigskip
\centerline{\sc Charles P. Boyer~~ Krzysztof Galicki~~ Michael Nakamaye}
\footnote{}{\ninerm During the preparation of this work the first two authors 
were partially supported by NSF grant DMS-9970904, and third author by NSF
grant DMS-0070190. 1991 Mathematics Subject Classification: 53C25,53C12,14E30}
\bigskip

\centerline{\vbox{\hsize = 5.85truein
\baselineskip = 12.5truept
\eightrm
\noindent {\bsc Abstract:}
We show that $\scriptstyle{\#9(S^2\times S^3)}$ admits an 8-dimensional
complex family of inequivalent non-regular \Se structures. These are the first
known Einstein metrics on this 5-manifold. In particular, the bound
$\scriptstyle{b_2(M)\leq8}$  which holds for any regular \Se $\scriptstyle{M}$
does not apply to the non-regular case. We also discuss the failure of the
Hitchin-Thorpe inequality in the case of 4-orbifolds and describe the orbifold
version.}} \tenrm

\vskip 1in
\baselineskip = 10 truept
\centerline{\bf Introduction}  
\bigskip
Recently Demailly and K\'ollar have developed some new techniques to study the
existence of K\"ahler-Einstein metrics on compact Fano orbifolds [DK].
Johnson and K\'ollar applied these techniques to study K\"ahler-Einstein
metrics on certain log del Pezzo surfaces in weighted projective
3-spaces [JK1] as well as anti-canonically embedded orbifold Fano
3-folds in weighted projective 4-spaces [JK2]. In [BG3, BGN1] we have extended
some of the results of [JK1] to the case of higher index and have studied
their implications in the realm of \Se metrics on simply connected
smooth 5-manifolds. These arise as links of isolated hypersurface singularities
given by quasi-homogeneous polynomials in $\bbc^4$ [BG3].  In [BGN1] we showed
that there are many families of non-regular \Se structures on $k$-fold
connected sums of $S^2\times S^3$ for $k=1,2,3,4,5,6,7$. When $k=8$ we have not
found any non-regular \Se structures. However, $8\#(S^2\times S^3)$ viewed as
a circle bundle over a blow-up of $\bbc\bbp^2$ at eight points comes with
an eight complex parameter family of 
regular \Se structures. It follows from the well-known classification of
smooth del Pezzo surfaces that for regular \Se 5-manifolds $\cals$ we must
have   $b_2(\cals)\leq8$. Moreover, it follows from the Hitchin-Thorpe
inequality that $\bbc\bbp^2$ blown-up at 9 or more points in general position
is excluded from admitting any Einstein metric whatsoever. Neither of these
bounds hold in the orbifold case. For the Hitchin-Thorpe type result the
topological Euler number and Hirzebruch signature must be replaced by their
corresponding orbifold analogue. Hence, no such topological bound holds in the
orbifold case. Nevertheless, all the non-regular examples found in
[BGN1] still satisfy this bound. The main purpose of this note is to show
that the bound does not hold for non-regular \Se 5-manifolds.
We achieve this by a careful analysis of the question of the existence
of K\"ahler-Einstein metric for one of the log del Pezzo surface $\calz_{16}
\subset\bbc\bbp^2(1,3,5,8)$ of
index 1 and degree 16 found by Johnson and K\'ollar [JK1]. Johnson and K\'ollar
leave the question of existence of K\"ahler-Einstein metrics on
$\calz_{16}$ open. Using well-known techniques of Milnor and Orlik [MO], 
one can easily compute  the characteristic polynomial of the associated link 
$L_{16}\subset \bbc^4$ and from this that the second Betti number
$b_2(L_{16})=9$. As a consequence of this, a general
argument concerning the absence of torsion in $H_2$ [BG3], and a
classification theorem of Smale [Sm], the link $L_{16}$, which is naturally a
$V$-bundle over $L_{16}\rightarrow\calz_{16},$ 
must be diffeomorphic to $9\#(S^2\times S^3)$. In section 3,
we show that $\calz$ does satisfy the sufficient conditions of
[DK] and does admit a K\"ahler-Einstein metric. This is the key to establishing 
the following

\noindent{\sc Theorem} A: \tensl Let $\cals=\#9(S^2\times S^3)$. 
Then $\cals$ admits an 8-dimensional family of inequivalent
non-regular Sasakian-Einstein structures. Furthermore, the Einstein metrics
in the family are inequivalent as Riemannian metrics.  \tenrm

It follows that the metric cone on $\cals=\#9(S^2\times S^3)$ is
is a singular Calabi-Yau 3-fold [BG1-2]. Such cones arise naturally in the
context of the supersymmetric string theory (see, for example, [Y]).
In view of Theorem A and Smale's classification result [Sm] 
one can ask the following important questions:

\noindent{\sc Problem} B: \tensl Suppose $M$ is a compact simply-connected
\Se 5-manifold. Can $b_2(M)$ be arbitrarily large? And, if so,
is there a \Se structure on any connected sum of $k$ copies of
$S^2\times S^3$?\tenrm

We are unable to answer this question here; however,
in a forthcoming note [BGN2] the authors show that there is a positive Sasakian
structure on $k\#(S^2\times S^3)$ for $k\geq1.$

\bigskip
\baselineskip = 10 truept
\centerline{\bf \sea. Sasakian Structures on Links of Isolated Hypersurface
Singularities}

Although the purpose of this not is to describe a solitary example
of the \Se geometry of $9\#(S^2\times S^3)$ we shall begin with
a very brief summary of the Sasakian and \Se geometry of links of isolated
hypersurface singularities defined by weighted homogeneous polynomials.
For more details we refer the reader to [BG3, BGN1].
Consider the affine space $\bbc^{n+1}$ together with a weighted
$\bbc^*$-action $\bbc^*_\bfw$ given by $(z_0,\ldots,z_n)\mapsto
(\grl^{w_0}z_0,\ldots,\grl^{w_n}z_n),$ where the {\it weights} $w_j$ are
positive integers. It is convenient to view the weights as the components of a
vector $\bfw\in (\bbz^+)^{n+1},$ and we shall assume that
$\gcd(w_0,\ldots,w_n)=1.$ Let $f$ be a quasi-homogeneous polynomial, that is
$f\in \bbc[z_0,\ldots,z_n]$ and satisfies
$$f(\grl^{w_0}z_0,\ldots,\grl^{w_n}z_n)=\grl^df(z_0,\ldots,z_n),
\leqno{\sea.1}$$
where $d\in \bbz^+$ is the degree of $f.$ We are interested in the {\it
weighted affine cone} $C_f$ defined by
the equation $f(z_0,\ldots,z_n)=0.$ We shall assume that the origin in
$\bbc^{n+1}$ is an isolated singularity, in fact the only singularity, of
$f.$ Then the link $L_f$ defined by 
$$L_f= C_f\cap S^{2n+1}, \leqno{\sea.2}$$
where 
$$S^{2n+1}=\{(z_0,\ldots,z_n)\in \bbc^{n+1}|\sum_{j=0}^n|z_j|^2=1\}$$
is the unit sphere in $\bbc^{n+1},$ is a smooth manifold of dimension $2n-1.$ 
Furthermore, it is well-known [Mil] that the link $L_f$ is $(n-2)$-connected.

On $S^{2n+1}$ there is a well-known  ``weighted'' Sasakian structure  
$\cals_\bfw=(\xi_\bfw,\eta_\bfw,\Phi_\bfw,g_\bfw)$ where the vector field
$\xi_\bfw$ is the infinitesimal generator of the circle subgroup $S^1_\bfw
\subset \bbc^*_\bfw.$ This Sasakian structure on $S^{2n+1}$ induces a
Sasakian structure, also denoted by
$\cals_\bfw$, on the link $L_f.$ (See
[YK, BG3, BGN1] for details. 
The quotient space $\calz_f$ of $S^{2n+1}$ by  $S^1_\bfw,$ or
equivalently the space of leaves of the characteristic foliation $\calf_\xi$
of $\cals_\bfw,$  is a compact K\"ahler orbifold which is a projective
algebraic variety embedded in the weighted projective $\bbp(\bfw)=
\bbp(w_0,w_1,\cdots,w_n),$ in such a way that there is a commutative diagram 
$$\matrix{L_f &\ra{2.5}& S^{2n+1}_\bfw&\cr
  \decdnar{\pi}&&\decdnar{} &\cr
   \calz_f &\ra{2.5} &\bbp(\bfw),&\cr}$$
where the horizontal arrows are Sasakian and K\"ahlerian embeddings,
respectively, and the vertical arrows are principal $S^1$ V-bundles and
orbifold Riemannian submersions.   

As with K\"ahler structures there are many Sasakian structures on a
given Sasakian manifold. In fact there are many Sasakian structures 
which have $\xi$ as its characteristic vector field. Such deformations of a
given Sasakian structure $\cals=(\xi,\eta,\Phi,g)$ are obtained by adding to
$\eta$ a  continuous one parameter family of basic1-forms $\grz_t.$          
We require that the 1-form $\eta_t=\eta +\grz_t$  satisfy the conditions  
$$\eta_0=\eta, \qquad \grz_0=0,\qquad \eta_t\wedge (d\eta_t)^n\neq 0~~
\forall~~ t\in [0,1]. \leqno{\sea.5}$$    
Since $\grz_t$ is basic $\xi$ is the Reeb (characteristic) vector field
associated to $\eta_t$ for all $t.$ Now let us define 
$$\eqalign{\Phi_t&=\Phi -\xi\otimes \grz_t\circ \Phi, \cr
           g_t&=g+d\grz_t\circ (\Phi\otimes \hbox{id}) + \grz_t\otimes\eta
+\eta\otimes \grz_t +\grz_t\otimes \grz_t.} \leqno{\sea.6}$$
Then $\cals_t=(\xi,\eta_t,\Phi_t,g_t)$ is a Sasakian structure for all $t\in
[0,1]$ that has the same underlying contact structure and the same
characteristic foliation. In general these structures are inequivalent and
the moduli space of Sasakian structures having the same characteristic vector
field is infinite dimensional. 

Suppose now we have  a link $L_f$ with a given
Sasakian structure $(\xi,\eta,\Phi,g)$. When can we find a 1-form $\grz$ such
that the deformed structure $(\xi,\eta+\grz,\Phi',g')$ is \Se? This is a
Sasakian version of the Calabi problem and its solution is equivalent to
solving the corresponding Calabi problem on the space of leaves $\calz_f$. 
Since a \Se manifold necessarily has positive Ricci tensor, its
Sasakian structure is necessarily positive. This also implies that the
K\"ahler structure on $\calz_f$ be positive, i.e. $c_1(\calz_f)$ can be
represented by a positive definite $(1,1)$ form. In this case there are 
well-known obstructions to solving to solving the Calabi problem.  These 
obstructions for finding a solution to the Monge-Ampere equations involves
the non-triviality of certain {\it multiplier ideal sheaves} [Na, DK] associated
with effective canonical $\bbq$-divisors on the space of leaves $\calz_f.$
Consequently, if one can show that these  multiplier ideal sheaves coincide
with the full structure sheaf, one obtains the existence of a positive
K\"ahler-Einstein metric on $\calz_f$ and hence, a \Se metric on $L_f.$ 
Sufficient conditions that guarantee this come from Mori theory and can be
phrased as follows [JK1]: 
\medskip
\noindent{\sc Sufficient Conditions for a K-E metric on $\calz_f$ and S-E
metric on $L_f$} \sea.7: 

\item{I.} For some $\gre > 0$ and every effective $\bbq$-divisor $D$  on
$\calz_f$ numerically equivalent  to $-K_{\calz_f},$  the pair 
$(\calz,{n+\gre\over n+1}D)$ is Kawamata log terminal (klt).

Now conditions on the weights that guarantee that the hypersurface $C_f\subset
\bbc^{n+1}$ have only an isolated singularity at the origin are well-known
[Fle, JK1].  These conditions, known as {\it quasi-smoothness} conditions
guarantee that $\calz_f$ is smooth in the orbifold sense, that is, at a vertex
 $P_i\in\bbp({\bf w})$ the preimage of $\calz_f$ in the orbifold chart of 
$\bbp({\bf w})$ is smooth.  It is easy to see that one can formulate all these
conditions as follows [Fle, JK1]:
\bigskip
\noindent{\sc Quasi-Smoothness Conditions} \sea.8:
\bigskip
\+II.& For each $i=0,\cdots,3$ there is a $j$ and a monomial
$z_i^{m_i}z_j\in \calo(d).$\cr
\+&Here $j=i$ is possible.\cr
\medskip
\+III.&  If $\gcd(w_i,w_j)>1$ then there is a
monomial $z_i^{b_i}z_j^{b_j}\in
\calo(d).$\cr
\medskip
\+IV.& For every $i,j$ either
there is a monomial $z_i^{b_i}z_j^{b_j}\in
\calo(d),$\cr
\+&or there are monomials $z_i^{c_i}z_j^{c_j}z_k$ and
$z_i^{d_i}z_j^{d_j}z_l\in \calo(d)$ with $\{k,l\}\neq \{i,j\}.$\cr

In condition I the $i=j$ case corresponds to the
case when $\calz_f$ does not pass through the point $P_i$. The second
condition is equivalent to $\calz_f$ not containing any of the singular
lines in $\bbp({\bf w})$. If $\calz_f$ contains a
coordinate axis (say $z_j=z_j=0$) then  the condition III
forces $\calz_f$ to be smooth
along it, except possibly at the vertices.

There is another condition apart from quasi-smoothness that assures us
that the adjunction theory behaves correctly, and that $\bbp({\bf w})$
does not have any orbifold singularities of codimension 1. It is [Dol, Fle]
\medskip
\noindent{\sc Well-formedness Condition} \sea.9
\bigskip
\+V.& For any triple $i,j,k\neq,$ we have $\gcd(w_i,w_j,w_k)=1.$ \cr
\medskip

Condition IV
guarantees that the canonical V-bundle $K_\calz$ is
determined in terms of the degree and index by
$$K_\calz \simeq \calo(-I)=\calo(d-|\bfw|),\leqno{\sea.10}$$
where $|\bfw|=\sum_iw_i.$

\bigskip
\baselineskip = 10 truept
\centerline{\bf \sec. Sufficient Conditions for K\"ahler-Einstein Metrics}
\bigskip

In this section we consider hypersurfaces $\calz_{16}$ of degree 16 and weights
$(1,3,5,8).$ As noted in [BGN1] there are 20 monomials in
$H^0({\bbp}(1,3,5,8),\calo(16)).$ To begin we consider the hypersurface
$\calz_{16}$ given by the zero set of  
$$
f(z_0,z_1,z_2,z_3)=z_0^{16} + z_1^5z_0 + z_2^3z_0 + z_1z_2z_3 + z_3^2
$$
in the weighted projective space ${\bbp}(1,3,5,8)$.  Thus $\calz_{16}$ is a 
hypersurface of degree 16 and index 1.  This hypersurface 
passes through
only two singular points of $\bbp(1,3,5,8)$, namely $P_1 = (0,0,1,0)$ and
$P_2 = (0,1,0,0)$.  We will show

\noindent{\sc Theorem} \sec.1 \tensl
For any effective divisor 
$$
D \in \left|\calo_{\calz_{16}}\left(-{11\over 16}K_{\calz_{16}}\right)\right|
$$
the pair $(\calz_{16},D)$ is klt.
\tenrm

\noindent
{\bf Proof:}  Suppose $D \in |-K_{\calz_{16}}|$.
If $Q \neq P_1,P_2$ then we see, since the 
linear series $\calo_{\calz_{16}}(3)$ has only isolated base points, that 
$$
{\rm mult}_Q(D) \leq D \cdot \calo_{\calz_{16}}(3) = 
{3 \cdot 16\over 3 \cdot 5
\cdot 8} = {2\over 5} < 1
$$
and hence $D$ is klt at $Q$ (see [BGN1] Lemma 2.11).  

We next consider the point $P_1$ which is the most difficult.
Let $\pi:(\bbc^2,0) \rightarrow (\calz_{16},P_1)$ 
be a local cover of the index 5 quotient 
singularity at $P_1$.  Intersecting with a general member of 
$\pi^\ast \calo_{\calz_{16}}(3)$, 
all of which pass through the point $P_1$, we see that 
$$
{\rm mult}_0(\pi^\ast D) \leq {5 \cdot 3 \cdot 16\over 3 \cdot 5 \cdot 8} = 2.
$$
Thus we may apply Shokurov's inversion of adjunction ([KM] Theorem
5.50) which says that the
pair $(\calz_{16},D)$ is log--canonical
at $P_1$ provided the pair $p_\ast^{-1}\pi^\ast D$ is a sum
of points with multiplicity at most one: here $p: S \rightarrow \bbc^2$ is the 
blow-up at the origin.  

In order to analyze $p_\ast^{-1}(\pi^\ast D)$ 
we first look at the natural candidate
for ``worst case scenario,'' namely $D = Z(z_0)$.  This is a curve 
$C \subset \bbp(3,5,8)$ with two components, namely $C_1 = Z(z_3)$ and
$C_2 = Z(z_3 + z_1z_2)$.  Let $D_1 = \pi^\ast C_1$ and $D_2 = \pi^\ast C_2$.
We claim that $D_1$ and $D_2$ are both smooth at $0$.  To see this, consider
$F = Z(z_1)$.  Then $F \cap C$ is proper and hence
$$
\pi^\ast F \cdot \pi^\ast C \geq 2
$$
with equality holding if and only if $\pi^\ast F$ 
meets $D_1$ and $D_2$ transversally
at $0$ (and at no other point).  It is easy to check that $F$ and $C$ meet
only at $P_1$ and hence $\pi^\ast F $ meets $D_1$ and $D_2$ only at $0$.
We compute the intersection multiplicity
$$
i_0(\pi^\ast F \cdot \pi^\ast C) = {5 \cdot 3 \cdot 16\over3 \cdot 5 \cdot
8} = 2.
$$
Thus $\pi^\ast F$ meets $D_1$ and $D_2$ transversally at $0$ and so $D_1$ and
$D_2$ are both smooth at $0$.

We now show that $D_1$ and $D_2$ meet transversally at $0$.  We already know
that $D_1$ meets $\pi^\ast F$ transversally.  Let $f_1$ be a local defining
equation of $D_1$ near $0$ and $f_2$ a local defining equation of $D_2$.  
Also, let $m_0$ denote the maximal ideal of the
local ring of $\bbc^2$ at $0$.  Then $f_1$ and $f_2$
generate $m_0$ if and only if $f_1$ and $f_1 - f_2$ generate $m_0$.  
But $Z(f_1) = \pi^\ast Z(z_3)$ and $Z(f_2) = \pi^\ast Z(z_1z_2 + z_3)$.  Thus
$Z(f_1 - f_2) = \pi^\ast Z(z_1z_2)$.  Since $0 \not\in \pi^\ast 
Z(z_2)$, we see that 
$D_1$ and $D_2$ meet transversally at $0$ if and only if $D_1$ and $\pi^\ast
(Z(z_1)) = \pi^\ast F$ meet transversally at $0$.  Thus, using the previous 
paragraph, we have established that $D_1$ and $D_2$ are both smooth at $0$
and meet transversally, i.e. the divisor $D = Z(z_0)$ is klt at $P_1$.

Now, following [JK1] and [BGN1], 
we deal with
an arbitrary (effective) choice of $D \in |-K_{\calz_{16}}|$.  Write
$$
D = aC_1 + bC_2 + D^\prime,\leqno{\sec.2}$$
where $D^\prime$ meets $C_1$ and $C_2$ properly.  We compute intersection
numbers:
$$\leqalignno{
\calo_{\calz_{16}}(1)\cdot(C_1+C_2) &= {16\over3\cdot 5\cdot 8} =
                      {2\over15},&\sec.3\cr 
\calo_{\calz_{16}}(1) \cdot C_1 &= {8\over 3 \cdot 5 \cdot 8} = {1\over 15},& 
\sec.4\cr
\calo_{\calz_{16}}(1) \cdot C_2 &={1\over15},& \sec.5\cr
C_1 \cdot C_2 &={1\over5},& \sec.6\cr 
C_1^2 = C_2^2 &= {-2\over15},& \sec.7\cr
C_1 \cdot D^\prime &\leq {2-a\over 15}.& \sec.8\cr}$$ 

The first inequality \sec.3 follows since ${z_0 = 0}$ has $C_1 + C_2$
as divisor of zeroes; \sec.4 follows since $C_1$ is the intersection of
the two hypersurfaces $z_0 = 0$ and $z_3 = 0$.  The third formula \sec.5
follows from \sec.3 and \sec.4.  Since $\pi^\ast C_1$
and $\pi^\ast C_2$ meet transversally at $0$ and nowhere else, this implies
\sec.6.  The next equation, \sec.7,
 follows from \sec.4, \sec.5,
and \sec.6,
using the fact that $\calo_{\calz_{16}}(1) \cdot \calz_{16}$ 
is represented by $C_1 + C_2$.  
Finally, \sec.8 follows since
$$
C_1 \cdot D^\prime \leq D \cdot D^\prime \leq D \cdot D - a D \cdot C_1.
$$

Intersecting \sec.2 with $C_1$ and using \sec.4, \sec.6,
and \sec.7 gives
$$\leqalignno{
{1\over15} &=C_1 \cdot \calo_{\calz_{16}}(1) \cr
             &= C_1(aC_1+bC_2+D^\prime) \cr
             &={-2a\over15} + {b\over5} + C_1 \cdot D^\prime.&\sec.9\cr}$$

Applying \sec.8 and using the fact that $a + b \leq 2$ (since
${\rm mult}_0(\pi^\ast D) \leq 2$), \sec.9 implies
$$
{2a\over15} \leq {1-a\over 15} + {2-a\over5}.\leqno{\sec.10}$$

Rearranging terms in \sec.10 gives
$$a \leq {7\over6}.\leqno{\sec.11}$$
An identical argument show that $b \leq {7\over6}$.

It remains to deal with the final term $D^\prime$.  Since
$D \cdot D^\prime \leq {2-a\over15}$, it follows that
$${\rm mult}_0(\pi^\ast D^\prime) \leq {2-a\over3}.\leqno{\sec.12}$$
Thus, in the tangent direction of $D_1$ at $0$ one obtains a divisor with
coefficient at most
$$
a + {2-a\over3} = {2a+2\over3} \leq {13\over9}.
$$
An identical bond holds with $b$ in place of $a$ and so we conclude that
$({\calz_{16}},{11\over16}D)$ is klt at $P_1$.
\medskip
We now turn to the point $P_2$ which will be simple as we have already
performed the relevant computations.  Let $\pi: (\bbc^2,0) \rightarrow
\calz_{16}$ be a local cover of the quotient singularity at $P_2$.  
Shokurov's 
inversion of adjunction still applies, as one sees by intersecting with 
a general member of $|\calo_{\calz_{16}}(5)|$.  We again let 
$Z(z_0) = C_1 \cup C_2$.
Taking $F = Z(z_2)$ and arguing
exactly as above establishes that $\pi^\ast(C_1)$ and $\pi^\ast(C_2)$ meet
transversally at $0$.  

If  $D \in |-K_{\calz_{16}}|$ is an effective divisor then as above we write
$$
D = aC_1 + bC_2 + D^\prime
$$
where $a, b \leq {7\over6}$.  The bound \sec.12 still holds for
${\rm mult}_0(\pi^\ast(D^\prime))$ (in fact, one can do better since the 
ramification index of $\pi$ is only 3) and so exactly as before we find that 
$\left({\calz_{16}},{11\over16}D\right)$ is klt at $P_2$.
\medskip
\medskip
Note that in the proof of Theorem \sec.1 the same argument works for
a much more general hypersurface $Y$.  In particular, any of the numerous
monomials with a $z_0$ term can be added without changing the proof.  The
only monomial which presents a problem is $z_1^2z_2^2$.  Consider the 
hypersurface $Y$ defined by 
$$
f(z_0,z_1,z_2,z_3) = z_0^{16} + z_1^5z_0 + z_2^3z_0 + z_3^2 + z_1^2z_2^2.
$$
The divisor $z_0 = 0$ again splits into two irreducible components,
$C_1  = Z(z_1z_2 + iz_3)$ and $C_2 = Z(z_1z_2 - i z_3)$. Arguing precisely
as above shows that the divisor $z_1 = 0$  meets $C_1$ and $C_2$ transversally
on the local cover $\pi: (\bbc^2,0) \rightarrow (Y,P_1)$.  
Hence $\pi^\ast(C_1^\prime)$ 
meets $\pi^\ast(C_2)$ transversally where $C_1^\prime = Z(z_1z_2)$.  Adding
the equations for $C_1$ and $C_2$ then establishes that 
$\pi^\ast(C_1)$ and $\pi^\ast(C_2)$ meet transversally. An identical
argument applies for the point $P_2$ as above.  
Thus the inversion of adjunction formula applies and the proof of Theorem 
\sec.1 can be copied to establish that $\left(Y,{11\over16}D\right)$ 
is klt for
any choice of effective $D \in |-K_Y|$.  
This argument always applies provided the quadratic form $as^2 + bst + ct^2$
has distinct roots where $a$ is the coefficient of $z_3^2$, $b$ is the 
coefficient of $z_1z_2z_3$ and $c$ is the coefficient of $(z_1z_2)^2$.
This then guarantees the existence of K\"ahler-Einstein metrics on the
hypersurfaces $\calz_{16}$ which depend on 20 complex parameters. We thus obtain
induced \Se metrics on the link $L_{16}$ depending on 20 complex parameters.
However, these parameters are not all effective, and in section \see~  we
discuss the moduli problem for these hypersurfaces and links. 

\bigskip
\baselineskip = 10 truept
\centerline{\bf \sed. The Topology of the Link $L_{16}$} 
\bigskip

The main purpose of this section is to prove

\noindent{\sc Theorem} \sed.1: \tensl  The \Se manifold $L_{16}$ is
diffeomorphic to $9\#(S^2\times S^3)$.\tenrm

In the process of proving this theorem we study the geometry of the link, and
in particular compute the  characteristic polynomial which is an important
link invariant that generalizes the Alexander polynomial of a knot. Let us 
recall the
well-known construction of Milnor [Mil] concerning isolated hypersurface
singularities: There is a fibration of $(S^{2n+1}-L_f)\ra{1.3} S^1$ whose
fiber $F$ is an open manifold that is homotopy equivalent to a bouquet of
n-spheres $S^n\vee S^n\cdots \vee S^n.$ The {\it Milnor number} $\mu$ of $L_f$
is the number of $S^n$'s in the bouquet. It is an invariant of the link which
can be calculated explicitly in terms of the degree $d$ and weights
$(w_0,\ldots,w_n)$ by the formula [MO]                                    
$$\mu =\mu(L_f)=\prod_{i=0}^n\bigl({d\over w_i}-1\bigr).\leqno{\sed.1}$$  
For our link $L_{16}$ with weights $\bfw=(1,3,5,8),$  one
immediately obtains

\noindent{\sc Proposition} \sed.2: \tensl The Milnor number of the 
simply connected \Se 5-manifold $L_{16}$ is
$\mu(L_{16})=143.$ \tenrm

The closure $\bar{F}$ of $F$ has the same homotopy type as $F$ and is a compact
manifold whose boundary is precisely the link $L_f.$ So the reduced homology of
$F$ and $\bar{F}$ is only non-zero in dimension $n$ and $H_n(F,\bbz)\approx
\bbz^{143}.$ Using the Wang sequence of the Milnor fibration together with
Alexander-Poincare duality gives the exact sequence [Mil]
$$0\ra{1.5} H_n(L_f,\bbz)\ra{1.5} H_n(F,\bbz)
\fract{\bbi -h_*}{\ra{1.5}} H_n(F,\bbz) \ra{1.5} H_{n-1}(L_f,\bbz)\ra{1.5} 0,
\leqno{\sed.3}$$ 
where $h_*$ is the {\it monodromy} map (or characteristic
map) induced by the $S^1_\bfw$ action. From this we see that
$H_n(L_f,\bbz)=\ker(\bbi -h_*)$ is a free Abelian group, and $H_{n-1}(L_f,\bbz)
=\coker(\bbi -h_*)$ which in general has torsion, but whose free part equals 
$\ker(\bbi -h_*).$ There is a well-known algorithm due to Milnor and Orlik [MO]
for computing the free part of $H_{n-1}(L_f,\bbz)$ in terms of the
characteristic polynomial $\grD(t)=\det(t\bbi -h_*),$ but finding the torsion
is, in general, much more difficult. However, in the case at hand, links of
dimension five, under a mild assumption a method of Randell [Ran] shows that
the torsion vanishes [BG3]. The Betti number $b_n(L_f)=b_{n-1}(L_f)$ equals the
number of factors of $(t-1)$ in $\grD(t).$ 

We now compute the characteristic polynomials $\grD(t)$ for our example.

\noindent{Proposition} \sed.4: \tensl The characteristic polynomial
$\grD(t)$ of the link
$L_{16}$ is given by 
$$\grD(t)= (t-1)^9(t+1)^8(t^2+1)^9(t^4+1)^9(t^8+1)^9.$$
Hence, the second Betti number is $b_2(L_{16})=9.$\tenrm

\noindent{\sc Proof}: The Milnor and Orlik [MO] algorithm for computing the
characteristic polynomial of the monodromy operator for weighted homogeneous
polynomials is as follows: First associate to any monic polynomial $F$ with
roots $\gra_1,\ldots,\gra_k\in \bbc^*$ its divisor 
$$\div F= <\gra_1>+\cdots+<\gra_k>$$
as an element of the integral ring $\bbz[\bbc^*]$ and let $\grL_n= \div
(t^n-1).$  The rational weights $w'_i$ used in [MO] are related to our integer
weights $w_i$ by $w_i'={d\over w_i},$ and we write the $w'_i={u_i\over v_i}$ in
irreducible form. So we have $w'_0=16, w'_1={16\over 3}, w'_2={16\over 5},  
w'_3=2.$ Then  the divisor of the characteristic polynomial
is  
$$\div \grD(t)=({\grL_{16}\over 5}-1)({\grL_{16}\over
3}-1)(\grL_{16}-1)(\grL_2-1) \leqno{\sed.5}$$ 
which upon using the relations $\grL_a\grL_b=\gcd(a,b)\grL_{lcm(a,b)}$ reduces
to 
$$\div \grD(t)=9\grL_{16}-\grL_2+1.$$
The characteristic polynomial $\grD(t)$ is then determined from its divisor by
$$\grD(t)=(t-1)\prod(t^j-1)^{a_j}, \leqno{\sed.6}$$
where $\div \grD(t)=1+\sum_ja_j\grL_j.$ This gives the result.
\hfill\za

Now the proof of Theorem \sed.1 follows by the Smale Classification Theorem
[Sm], since $L_{16}$ has no torsion by Lemma 5.8 of [BG3], and any \Se
manifold is spin.

\bigskip
\baselineskip = 10 truept
\centerline{\bf \tvo. Topological versus Orbifold Invariants} 
\bigskip

It is well-known that dimensions three and four are the only dimensions where
there are known topological obstructions to the existence of Einstein metrics
on smooth manifolds. The first example of this in dimension four is due to
Marcel Berger [Be] who noticed that a compact four dimensional Einstein
manifold must have non-negative Euler characteristic. Later Thorpe [Be] and
then Hitchin [Hit] noticed a much sharper inequality obstructing Einstein
metrics.  

In this section we  discuss both topological and orbifold invariants of
compact complex four dimensional orbifolds with isolated orbifold
singularities, and the failure of the Hitchin-Thorpe (cf. [Be, Hit, LeB])
inequality in the singular case. Of course, we certainly don't  expect this to
hold as is in the orbifold case. There are orbifold correction terms as
indicated by Satake [Sat] for the Euler characteristic. We note that for any
log del Pezzo surface $\calz$ that both the topological Euler characteristic
and Hirzebruch signature are determined in terms of the second Betti number
$b_2(\calz)$ by
$$\chi_{top}(\calz)=2+b_2(\calz), \qquad \grt_{top}(\calz)=2-b_2(\calz).
\leqno{\tvo.1}$$
This gives
$$2\chi_{top}(\calz)+3\tau_{top}(\calz)=10-b_2(\calz). \leqno{\tvo.2}$$
In the smooth case this is positive for del Pezzo surfaces and any smooth rational surface with $b_2\geq 10$ cannot admit any Einstein metric. However, for
our singular $\calz_{16}$ we have  $$b_2(\calz_{16})=10, \qquad
\chi_{top}(\calz_{16})=12, \qquad \tau_{top}(\calz_{16})=-8. \leqno{\tvo.3}$$  
Thus, we have
$$2\chi_{top}(\calz_{16})+3\tau_{top}(\calz_{16})=0\leqno{\tvo.4}$$
even though $\calz_{16}$ has a positive K\"ahler-Einstein metric. Thus, the
Hitchin-Thorpe inequality clearly fails in the singular category. 

The point is, of course, that for orbifolds both the Gauss-Bonnet [Sat] and
signature theorems [Kaw] relate integrals of curvature invariants to rational
numbers that are orbifold invariants. In particular, the orbifold Euler
characteristic $\chi_{orb}(\calz)$ for a compact 4-orbifold $\calz$ with
at most isolated orbifold singularities  is given by [Sat]
$$\chi_{orb}(\calz)={1\over
8\pi^2}\int_{\calz}\Bigl(|W_+|^2 +|W_-|^2+{s^2\over 24}-{|\Ric^0|^2\over
2}\Bigr)\leqno{\tvo.5}$$ 
where $W_\pm,s,\Ric^0$ are the self-dual (anti-self-dual) pieces of the Weyl
tensor, the scalar curvature, and traceless Ricci curvature, respectively. In
the case that $\calz$ is a compact complex surface $\chi_{orb}(\calz)$
equals [Bla] the top Chern class $c_2(\calz).$ Moreover,
$\chi_{orb}$ is related to the topological Euler characteristic
$\chi_{top}(\calz)$ by [Bla, Sat]  
$$\chi_{orb}(\calz)=\chi_{top}(\calz)-\sum_x\Bigl(1-{1\over
|\grG_x|}\Bigr) \leqno{\tvo.6}$$ 
where $|\grG_x|$ denotes the order of the local uniformizing
group $\grG_x$ at the singular point $x$ and the sum is taken over all
singular points. Notice that for orbifolds the analogue of Berger's result is
a bit stronger:

\noindent{\sc Proposition} \tvo.7: \tensl Let $\calz$ be a compact 
4-dimensional orbifold with precisely $N$ isolated orbifold singularities.
Suppose further that $\calz$ admits an orbifold Einstein metric, then
$$\chi_{top}(\calz)\geq \sum_x\bigl(1-{1\over |\grG_x|}\bigr)\geq {N\over
2}.$$ \tenrm

A similar analysis can be done for the signature. The Hirzebruch signature
theorem says that for smooth manifolds whose dimension is a multiple of four
the signature $\tau$ can be obtained by evaluating the L-function of the
tangent bundle on the fundamental class of the  manifold. In the orbifold case
this is replaced by Kawasaki's signature theorem [Kaw], which gives
$$\tau_{top}(\calz)= <L(T\calz),[\calz]> + \sum_x{\call_x\over |\grG_x|},\leqno{\tvo.8}$$
where again the sum is over all singular points, and $\call_x$ denotes the
evaluation of the residual L-class at $x.$ Here we have
$$<L(T\calz),[\calz]>= {1\over 12\pi^2}\int_\calz\Bigl(|W_+|^2-|W_-|^2\Bigr).
\leqno{\tvo.9}$$

We can view this expression as an ``orbifold signature'' which is a rational
number, and write $\grt_{orb}(\calz)= <L(T\calz),[\calz]>.$ Then we rewrite
equation \tvo.8 as 
$$\grt_{top}(\calz)= \grt_{orb}(\calz)+ \grt_{res}(\calz)$$
and the analogue for orbifolds of the Hitchin-Thorpe curvature integral becomes
$$2\chi_{orb}(\calz)\pm 3\tau_{orb}(\calz) ={1\over
4\pi^2}\int_{\calz}\Bigl(2|W_\pm|^2+{s^2\over 24}-{|\Ric^0|^2\over
2}\Bigr),\leqno{\tvo.10}$$ 
so we obtain

\noindent{\sc Proposition} \tvo.11: \tensl Let $\calz$ be a compact 
4-dimensional orbifold with only isolated orbifold singularities.
Suppose further that $\calz$ admits an orbifold Einstein metric, then
$$\chi_{orb}(\calz)\geq {3\over 2}|\tau_{orb}(\calz)|.$$ 
Thus, we have
$$\chi_{top}(\calz)\geq {3\over 2}|\tau_{orb}(\calz)|+ \sum_x\bigl(1-{1\over
|\grG_x|}\bigr)\geq {3\over 2}|\tau_{orb}(\calz)|+ {N\over 2}.$$
\tenrm

In general Proposition \tvo.11 contains much less information than the
Hitchin-Thorpe result, since it involves orbifold invariants and not
topological invariants. Nevertheless, it does give an improvement of
Proposition \tvo.7.  We also have

\noindent{\sc Proposition} \tvo.12: \tensl Let $\calz$ be a compact complex
4-orbifold. Then
$$c_1^2(\calz)=2\chi_{orb}(\calz)+3\tau_{orb}(\calz).$$ 
Thus, if $\calz$ admits an Einstein metric compatible with the orbifold
structure, then $c_1^2\geq 0.$
\tenrm

\noindent{\sc Proof}: The Chern number $c_1^2$ is obtained by evaluating the
square of the Ricci form $\rho$ on the fundamental class of the orbifold, and
the usual proof holds once we replace the topological invariants $\chi_{top}$
and $\tau_{top}$ by their orbifold counterparts.  The last statement then
follows immediately from Proposition \tvo.11. \hfill\za 

\noindent{\sc Remark} \tvo.13: If $\calz$ is a compact complex 4-orbifold
embedded as a hypersurface of degree $d$ in the weighted projective space
$\bbp(\bfw),$ then it is easy to see that 
$$c_1^2={d(|\bfw|-d)^2\over w_0w_1w_2w_3}.\leqno{\tvo.14}$$
Thus, we see that if a compact complex 4-orbifold $\calz$ has
$2\chi_{orb}+3\tau_{orb}=c_1^2<0$ it cannot admit a compatible Einstein metric
nor be embedded as a hypersurface in $\bbp(\bfw).$ We also mention that the
case that the canonical V-bundle is trivial, so $c_1^2=0$  and $\calz$ is
simply connected gives Reid's list [Fle] of 95 mostly singular K3 surfaces.  

Let us now evaluate the orbifold invariants for the case currently under
study, namely the singular complex surface $\calz_{16}.$ 
In this case it is easy to evaluate both $c_1^2$ and $\chi_{orb}$ directly
giving 
$$c_1^2(\calz_{16})={16\over 3\cdot 5\cdot 8}={2\over 15},\qquad
\chi_{orb}(\calz_{16})=12- {2\over 3}-{4\over 5}= 10+{8\over 15}.$$
We have not evaluated $\tau_{orb}(\calz_{16})$ directly, but \tvo.12 gives 
$\tau_{orb}(\calz_{16})=-7 +{1\over 45}.$ 

\bigskip
\baselineskip = 10 truept
\centerline{\bf \see. The Moduli of Sasakian-Einstein Structures on $L_{16}$}
\bigskip

To analyze the moduli problem in the non-regular case, we
begin by describing the group of complex automorphism
$\gG_\bfw$ of the weighted
projective 3-space $\bbp(\bfw)$. We shall assume that
$\bbp(\bfw)$ is well-formed. Recall that
$\bbp(\bfw)$ can be defined as a scheme ${\rm Proj}(S(\bfw))$,
where
$$S(\bfw)=\bigoplus_dS^d(\bfw)=\bbc[z_0,z_1,z_2,z_3].$$
The ring of polynomials $\bbc[z_0,z_1,z_2,z_3]$ is graded with grading
defined by the weights $\bfw=(w_1,w_1,w_2,w_3)$. As a projective variety
we can embed $\bbp(\bfw)\subset \bbc\bbp^N$ and then the
group  $\gG_\bfw$ is a subgroup of $PGL(N,\bbc)$. $\bbp(\bfw)$ is
a toric variety and we can describe $\gG_\bfw$ explicitly as follows:
Let $\bfw=(w_0,w_1,w_2,w_3)$ be ordered as before.
We consider the group $G(\bfw)$ of automorphisms of the graded ring $S(\bfw)$
defined on generators by
$$\varphi_\bfw\pmatrix{z_0\cr z_1\cr z_2\cr z_3\cr}=
\pmatrix{f_0^{(w_0)}(z_0,z_1,z_2,z_3)\cr
f_1^{(w_1)}(z_0,z_1,z_2,z_3)\cr f_2^{(w_2)}(z_0,z_1,z_2,z_3)\cr
f_3^{(w_3)}(z_0,z_1,z_2,z_3)\cr},\leqno{\see.1}$$
where $f_i^{(w_i)}(z_0,z_1,z_2,z_3)$ is an arbitrary weighted homogeneous
polynomial of degree $w_i$ in $(z_0,z_1,z_2,z_3)$. This is
a finite dimensional Lie group and
it is a subgroup of $GL(N,\bbc)$.
Projectivising, we get $\gG_\bfw=\bbp_\bbc(G(\bfw))$.

Note that when $\bfw=(1,1,1,1)$ then $G(\bfw)=GL(4,\bbc)$.  Other than this
case three weights are never the same if $\bbp(\bfw)$ is well-formed. If two
weights coincide then $G(\bfw)$ contains $GL(2,\bbc)$ as a subgroup. Finally,
when all weights are distinct we can write
$$\varphi_\bfw\pmatrix{z_0\cr z_1\cr z_2\cr z_3\cr}=
\pmatrix{a_0z_0\cr
a_1z_1+f_1^{(w_0)}(z_1)\cr a_2z_2+f_2^{(w_2)}(z_0,z_1)\cr
a_3z_3+f_3^{(w_3)}(z_0,z_1,z_2)\cr},\leqno{\see.2}$$
where $(a_0,a_1,a_2,a_3)\in(\bbc^*)^4$ and $f_i^{(w_i)},\ \ i=1,2,3$
are weighted homogeneous polynomials of degree $w_i$. The simplest situation
occurs when $f_1=f_2=f_3$ are forced to vanish. Then $\gG_\bfw=(\bbc^*)^3$ is
the smallest it can possibly be as $\bbp(\bfw)$ is toric.

Let $S_\bfw^d\subset S(\bfw)$ be the vector subspace spanned by all monomials
in $(z_0,z_1,z_2,z_3)$ of degree $d=|w|-I$, and
let $\hat{S}^d(\bfw)\subset S^d(\bfw)$ denote subset all quasi-smooth
elements.
Then we define $m_\bfw^d$ to be the dimension of the subspace
generated by $\hat{S}_\bfw^d.$ Now the automorphism group $G(\bfw)$ acts on
$S_\bfw^d$ leaving the subset $\hat{S}^d(\bfw)$ of quasi-smooth polynomials
invariant. Thus, for each log del Pezzo surface
we define the moduli space
$$\calm_\bfw^d=\hat{S}_\bfw^d/G(\bfw)=\bbp(\hat{S}_\bfw^d)/\gG_\bfw,
\leqno{\see.3}$$
with $n_\bfw^d={\rm dim}(\calm_\bfw^d).$ Now there is an injective map
$$\calm_\bfw^d\ra{2.0} \calm^\bbc(\calz_\bfw),\leqno{\see.4}$$
and each element in $\calm_\bfw^d$
corresponds to a unique homothety class of K\"ahler-Einstein metrics modulo
$\gG_\bfw$ and hence, to a unique \Se structure on the corresponding 5-manifold
$\cals_l$ modulo the group $\gG_\bfw$ acting as CR automorphisms. 

Let us now consider our $\calz_{16}$. It is easy to see that
$S^{16}(1,3,5,8)$ is isomorphic to $\bbc^{20}$ 
and it is spanned by the monomials
$z_3^2,z_1z_2z_3$, $z_1^2z_2^2$, $z_0z_2^3$, $z_0z_1^5$, $z_0^2z_1^2z_3$,
$z_0^2z_1^3z_2$, $z_0^3z_2z_3$, $z_0^3z_1z_2^2$, $z_0^4z_1^4$, $z_0^5z_1z_3$,
$z_0^5z_1^2z_2$, $z_0^6z_2^2$, $z_0^7z_1^3$, $z_0^8z_3$, $z_0^8z_1z_2$,
$z_0^{10}z_1^2$, $z_0^{11}z_2$, $z_0^{13}z_1$, $z_0^{16}.$ We take the
open submanifold $\hat S^{16}(1,3,5,8)\subset S^{16}(1,3,5,8)$.
This is acted on by the complex automorphism group, namely the group
$G(1,3,5,8)$ generated by
$$\varphi_\bfw\pmatrix{z_0\cr z_1\cr z_2\cr z_3\cr}=
\pmatrix{a_0^1z_0\cr
a_1^1z_1+a_1^2z_0^3\cr a_2^1z_2+a_2^2z_0^5+a_2^3z_0^2z_1\cr
a_3^1z_3+a_3^2z_0^8
+a_3^3z_1z_0^5+a_3^4z_1^2z_0^2+a_3^5z_2z_0^3+a_3^6z_2z_1\cr},\leqno{\see.5}$$
where $a_i^1\in\bbc^*$ and all other coefficients are in $\bbc.$
$G(1,3,5,8)$ is a 
12-dimensional complex Lie group acting on the open submanifold
$\hat S^{16}(1,3,5,8)\subset S^{16}(1,3,5,8)\approx\bbc^{20}.$ It follows 
that the quotient is an 8-dimensional complex manifold which by the
Bando-Mabuchi Theorem [BM] is the moduli space of positive K\"ahler-Einstein
metrics on the underlying compact orbifold $\calz_{16}$.

Theorem A now follows from these result and Proposition 7.14 of [BGN1].
%\vfil\eject
\bigskip
\medskip
\centerline{\bf Bibliography}
\medskip
\font\ninesl=cmsl9
\font\bsc=cmcsc10 at 10truept
\parskip=1.5truept
\baselineskip=11truept
\ninerm

\item{[Be]} {\bsc A. Besse}, {\ninesl Einstein manifolds},
Springer-Verlag, Berlin-New York, 1987.
\item{[BG1]} {\bsc C. P. Boyer and  K. Galicki}, {\ninesl On Sasakian-Einstein
Geometry}, Int. J. of Math. 11 (2000), 873-909.
\item{[BG2]} {\bsc C. P. Boyer and  K. Galicki}, {\ninesl 3-Sasakian
Manifolds}, to appear in Essays on Einstein Manifolds, Surveys in Differential
Geometry Vol V, International Press, 2000, C. LeBrun and M. Wang, Eds.
\item{[BG3]} {\bsc C. P. Boyer and  K. Galicki}, {\ninesl New Einstein Metrics
in Dimension Five},  math.DG/0003174, submitted for publication.
\item{[BGN1]} {\bsc C. P. Boyer, K. Galicki, and M. Nakamaye}, {\ninesl On the
Geometry of Sasakian-Einstein 5-Manifolds}, math.DG/0012047; submitted
for publication.
\item{[BGN2]} {\bsc C. P. Boyer, K. Galicki, and M. Nakamaye}, {\ninesl On 
Positive Sasakian Geometry}, in preparation.
\item{[Bla]} {\bsc R. Blache}, {\ninesl Chern classes and
Hirzebruch-Riemann-Roch theorem for coherent sheaves on complex projective
orbifolds with isolated singularities}, Math. Z. 222 (1996), 7-57.
\item{[BM]} {\bsc S. Bando and T. Mabuchi}, {\ninesl Uniqueness of Einstein
K\"ahler Metrics Modulo Connected Group Actions}, Adv. Stud. Pure Math. 10
(1987), 11-40.
\item{[DK]} {\bsc J.-P. Demailly and J. Koll\'ar}, {\ninesl Semi-continuity of
complex singularity exponents and K\"ahler-Einstein metrics on Fano
orbifolds}, preprint AG/9910118, to appear in Ann. Scient. Ec. Norm. Sup. Paris
\item{[Dol]} {\bsc I. Dolgachev}, {\ninesl Weighted projective varieties}, in
Proceedings, Group Actions and Vector Fields, Vancouver (1981) LNM 956, 34-71.
\item{[Fle]} {\bsc A.R. Fletcher}, {\ninesl Working with weighted complete
intersections}, Preprint MPI/89-95, revised version in  {\it Explicit
birational geometry of 3-folds},  A. Corti and M. Reid, eds.,
Cambridge Univ. Press, 2000,  pp 101-173.
\item{[Hit]} {\bsc N. Hitchin}, {\ninesl On compact four-dimensional Einstein
manifolds}, J. Diff. Geom. 9 (1974), 435-442.
\item{[JK1]} {\bsc J.M. Johnson and J. Koll\'ar}, {\ninesl K\"ahler-Einstein
metrics on log del Pezzo surfaces in weighted projective 3-space}, preprint
AG/0008129, to appear in Ann. Inst. Fourier.
\item{[JK2]} {\bsc J.M. Johnson and J. Koll\'ar},
{\ninesl Fano hypersurfaces in
weighted projective 4-spaces}, preprint AG/0008189, 
to appear in Experimental Math.
\item{[Kaw]} {\bsc Y. Kawamata}, {\ninesl  The cone of curves of algebraic
varieties}, Annals of Math 119, 1984, pp. 603-33.
\item{[KM]} {\bsc J. Koll\'ar, and S. Mori}, {\ninesl Birational Geometry of
Algebraic Varieties}, Cambridge University Press, 1998.
\item{[LeB]} {\bsc C. LeBrun}, {\ninesl Four-Dimensional Einstein Manifolds,
and Beyond}, Surveys in Differential Geometry VI:
{\it Essays on Einstein Manifolds};
A supplement to the Journal of Differential Geometry, pp. 247-285,
(eds. C. LeBrun, M. Wang); International Press, Cambridge (1999).
\item{[Mil]} {\bsc J. Milnor}, {\ninesl Singular Points of Complex
Hypersurfaces}, Ann. of Math. Stud. 61, Princeton Univ. Press, 1968.
\item{[MO]} {\bsc J. Milnor and P. Orlik}, {\ninesl Isolated singularities
defined by weighted homogeneous polynomials}, Topology 9 (1970), 385-393.
\item{[Na]} {\bsc A.M. Nadel}, {\ninesl Multiplier ideal sheaves
and existence of K\"ahler-Einstein metrics of positive scalar curvature}, Ann.
Math. 138 (1990), 549-596.
\item{[Ran]} {\bsc R.C. Randell} {\ninesl The homology of generalized
Brieskorn manifolds}, Topology 14 (1975), 347-355.
\item{[Sm]} {\bsc S. Smale}, {\ninesl On the structure of 5-manifolds},
Ann. Math. 75 (1962), 38-46.
\item{[Sat]} {\bsc I. Satake}, {\ninesl The Gauss-Bonnet theorem for  
$V$-manifolds}, J. Math. Soc. Japan V.9 No 4. (1957), 464-476. 
\item{[Y]} {\bsc S. -T. Yau}, {\ninesl Einstein manifolds with zero Ricci
curvature},
Surveys in Differential Geometry VI:
{\it Essays on Einstein Manifolds};
A supplement to the Journal of Differential Geometry, pp.1-14,
(eds. C. LeBrun, M. Wang); International Press, Cambridge (1999).
\item{[YK]} {\bsc K. Yano and M. Kon}, {\ninesl
Structures on manifolds}, Series in Pure Mathematics 3,
World Scientific Pub. Co., Singapore, 1984.
\medskip
\bigskip \line{ Department of Mathematics and Statistics
\hfil February 2001} \line{ University of New Mexico \hfil }
\line{ Albuquerque, NM 87131 \hfil } \line{ email: cboyer@math.unm.edu,
galicki@math.unm.edu, nakamaye@math.unm.edu\hfil} \line{ web pages:
http://www.math.unm.edu/$\tilde{\phantom{o}}$cboyer,
http://www.math.unm.edu/$\tilde{\phantom{o}}$galicki \hfil}
\end